\documentclass{article} 
\usepackage{amssymb, amsmath, amsthm,enumerate}
\theoremstyle{plain}
\newtheorem{theorem}{Theorem}
\newtheorem*{theorem*}{Theorem}
\newtheorem*{corollary*}{Corollary}

\newtheorem*{proposition*}{Proposition}
\newtheorem*{namedtheorem}{\theoremname}
\newcommand{\theoremname}{te}
\newenvironment{named}[1]{\renewcommand{\theoremname}{#1}
\begin{namedtheorem}}{\end{namedtheorem}}

\title{Alexey Vasilyevich Pogorelov,\\ the mathematician of an incredible power}
\author{A.A.Borisenko\\[2ex]
Kharkov Karazin National University,
\\e-mail: borisenk@univer.kharkov.ua }
\begin{document}
\maketitle
\begin{abstract}
Life and the mathematical legacy of the great mathematician A.V.~Pogorelov.
\footnote{This paper is the English translation of \cite{RU} }
\\[2ex]
\end{abstract}

\section*{Introduction.}

By the beginning of the 20th century, the local differential geometry of surfaces was very well-developed, mostly
by the means of the local analysis. In contrast, there was an apparent lack of methods and results in the global
geometry, ``geometry in the large", in which both geometry and analysis were equally helpless.
A typical example of such a question is the classical problem of the rigidity of a closed convex surface (ovaloid).
The contribution to this problem was made by brilliant mathematicians, such as Liebmann, Minkowski, Hilbert,
Weyl, Blaschke. A major progress was achieved by Cohn-Vossen at the beginning of the 1920-th. He proved that
isomeric $C^3$ ovaloids of positive Gauss curvature are congruent.
Meanwhile, by a classical result of Cauchy, isomeric closed convex polyhedrons are congruent. It seemed that
these two results are the particular cases of a general theorem of congruency
of any two (generally speaking, irregular) isomeric ovaloids. No approaches to the problem in such a generality
were visible. Similarly, it was not clear how to prove an isometric deformability of an ovaloid with a part removed,
and to estimate the deformability of an incomplete ovaloid. Under an assumption of a sufficient regularity, these
problems can be formulated in the language of nonlinear PDE's, but the theory of such equations was far from being
well developed at that time.

In many cases, analytic tools were specifically developed to study a particular geometric
problem. This is illustrated by the problem of the existence of a closed convex surface with a prescribed
analytic metric of positive Gauss curvature defined on a topological sphere. The general approach given by Weyl in
1916 was successfully completed in the 1930-th  by G.Levi who developed a very delicate techniques from the analytic
theory of the Monge-Amp\`{e}re equations. However, in the G.Levi's papers, analysis was developed separately from
geometry and
was applied to geometry as a ready-made tool. Many other papers just gave ad hoc methods for separate problems.
In fact, at that time, the fundamental problems of deformation of surfaces and many other problems of global
geometry remained unapproachable.

Another area of active research at the beginning of the 20th century was the theory of convex bodies in the Euclidean
space, including various geometric properties, the mixed volumes and inequalities. The cornerstone book
``Theory of convex bodies" by T.Bonnezen and T.Fenchel's  published in 1934 in German contained a complete
up-to-date account of the research in that area and an extensive bibliography. This book has not lost its value
till nowadays. In 2001, it was translated into Russian by V.A.Zalgaller. However, the study of convex surface was
beyond the scope of this book, and most probably, outside of the area of interest of researchers at that time.

\section {Intrinsic geometry of convex surfaces.}

Early papers of A.D.Aleksandrov were focused on the same questions. Starting with the classical results of
Minkowski, Alexander Danilovich Aleksandrov established new inequalities for the mixed volumes
of convex bodies. Curiously enough, forty years later, the algebraic analogues of his inequalities obtained as a
byproduct of his study, were successfully applied to the well-known Van-der-Waerden problem on the estimation of the
permanent, posed in 1926. Aleksandrov's inequalities for the mixed volumes also have interesting generalizations and
applications in algebraic geometry, in the theory of nonlinear elliptic equations and even in stochastic processes.

Simultaneously, A.D.Aleksandrov applied  the tools from the measure theory and functional analysis to the theory of
convex bodies. He introduced a functional space generated by the support functions and special measures on them,
the ``surface functions" and the related ``curvature functions". He proved that a convex body is determined uniquely,
up to translation, by the curvature function. This theorem includes, as the limiting cases, the theorems of
Christoffel and of Minkowski. In the course of proof, Alexander Danilovich introduced the concept of the generalized
differential equation in measures and of the corresponding generalized solution.

In 1941, A.D.Aleksandrov began the study of the intrinsic geometry of convex surfaces. He widened the class of
regular convex surfaces to the class of arbitrary convex surfaces (defined as domains on the boundary of a convex
body). The problems within that wider class required new techniques, beyond the Gaussian geometry of regular
surfaces.
It was necessary to understand the intrinsic geometry of an arbitrary convex surface (the properties
depending on the measurements on the surface itself) and to develop the tools for studying these properties, and then
to find the connection between the intrinsic and the extrinsic geometry of an arbitrary convex surface.
A.D.Aleksandrov constructed the intrinsic geometry of convex (and of arbitrary) surfaces starting from the general
concept of the metric space. Let $R$ be a metric space, a set such that for each pair of elements $X,Y \in R$ there
defined a number $\rho(X,Y)$, the \emph{distance}, satisfying the following axioms:
\begin{enumerate}
\item $ \rho (X, Y) \geqslant 0$ and $\rho (X, Y) =0 $ if and only if $X = Y.$
\item $\rho (X, Y)=\rho(Y,X)$.
\item $ \rho (X, Y) + \rho (Y, Z) \geqslant \rho (X, Z)$ (the triangle inequality).
\end{enumerate}
For example, the Euclidean space with the usual distance between the points is a metric space.
A \emph{curve} $ \gamma $ in a metric space $R$ is a continuous image of a segment $[a,b]$ considered together with
the continuous mapping $F:[a,b] \to R$. Just as in the for the Euclidean space, one can introduce the concept of the
length of a curve in a metric space by defining
$$
l _ {\gamma} = \sup \Sigma \, \rho (F (t _ {k-1}), \, F (t _ {k})), \ \ a=t_0 \leqslant t _ {1}
\leqslant t _ {2} \dots\leqslant t_n=b,
$$
where $ \rho (F (t _ {k-1}), \, F (t _ {k})) $ is the distance between $F (t_{k-1}) $ and $F (t_{k})$ in $R$, and the
supremum is taken over all finite partitions of the segment $ [a, b] $ by the $t_{k}$'s. The length so defined is
additive: if a curve $ \gamma $ is composed by the curves $ \gamma _ {1} $ and $ \gamma _ {2} $, then
$l(\gamma)=l(\gamma_1)+l(\gamma_2)$. Similar to the Euclidean space, in a metric space, the set of curves
of bounded length in a compact domain is compact, that is, any infinite sequence of such curves contains a converging
subsequence. Moreover, the length of the limiting curve is not greater than the lower limit of the lengths of curves
of the subsequence. Suppose that any two points $X, Y $ in $R $ can be connected by a rectifiable curve. Then
the intrinsic distance $\rho^{*}(X,Y)$ between $X, Y$ in $R$ can be defined as the infimum of the lengths of all the
rectifiable curves connecting $X$ and $Y$. It is easy to see that $\rho^*$ satisfies the axioms 1,2,3. The metric
$\rho^{*} $ is called the \emph{intrinsic} metric on $R$. If $\rho (X, Y) = \rho^{*}(X, Y)$ for any $X $ and $Y $,
then $R $ is called \emph{a space with an intrinsic metric} (the length space).
For a metric space $R$ to be isometric to a surface in the Euclidean space, its metric must necessarily be intrinsic.

Let $R $ be a manifold with an intrinsic metric. A curve $\gamma$ in $R$ is called \emph{a shortest path} if its length
is equal to the distance between the endpoints, hence being not greater than the length of any other curve joining the
same points. Each segment of a shortest path is also a shortest path. The limiting curve for a converging sequence of
shortest paths is again a shortest path.

In general, not every two points on a manifold can be connected by a shortest path, but every
point of a manifold has a neighborhood such that any two points from it can be joined by a shortest path. If a
manifold $R$ is metric complete, that is, if any closed bounded subset of it is compact,
then any two points of $R$ can be connected by a shortest path.
In  $R$, one can define a triangle, a polygon, a polygonal line, etc. in the usual way. The definition of the angle
between shortest paths is fundamental for the theory. The idea is to compare a triangle in $R$ to a triangle with the
same sides on the plane. Let $O \in R$ and let $OA$ and $OB$ be the shortest paths.
Choose arbitrary points $X \in OA$ and $Y \in OB$. Let $O'X'Y'$ be a triangle on the plane with the same sides as
$OXY$. Let $x=\rho (O,X), \; y=\rho (O, Y)$ and let $\alpha (x, y) $ be the angle at $O'$ of $O'X'Y'$.
The \emph{upper angle between the shortest paths $OX $ and $OY $ } is defined as the upper limit of $\alpha(x,y)$ when
$x,y \to 0$.

Using the notion of the angle, one can define the excess of a triangle as $ \alpha + \beta + \gamma-\pi $,
where $ \alpha, \beta, \gamma $ are the upper angles between the sides. At this point, one can already define
two-dimensional spaces of non-negative curvature. Namely, $R$ is called \emph{a space of non-negative curvature} if any
point of $R$ has a neighborhood such that excess of each triangle lying in it is non-negative. Each convex surface is
a space of non-negative curvature. The space $R$ has non-negative curvature if and only if
the function $\alpha (x,y) $ is non-increasing: if $x_1\leqslant x_2$ and $y_1\leqslant y_2 $, then
$\alpha (x_1, y_1) \geqslant \alpha (x_2, y_2) $. From the monotonicity of $\alpha$ it follows that the upper limit can
be replaced by the usual limit, which gives \emph{the angle between the shortest paths $OX $ and $OY$}.

The next step is to define the curvature. In the sense of A.D.Aleksandrov, the curvature is an additive function
of sets. The \emph{extrinsic curvature} of a set $M$ on a convex surface is the area of the spherical image
of $M$. The \emph{intrinsic curvature}, as an object of the intrinsic geometry of a convex surface, is
first defined for the three ``basic sets": 1) the curvature of an open triangle is its excess;
2) the curvature of an open shortest path is zero; 3) the curvature of a point equals $2\pi-\theta$, where $\theta$ is
the full angle around the point.
Then the curvature is (uniquely) defined by additivity for all Borel sets. A.D.Aleksandrov proved
that \emph{the curvature of any Borel set on a convex surface equals the area of its spherical image}, the
Gauss's Theorema Egregium for arbitrary convex surfaces. He established the fundamental connection
between the intrinsic geometry and the extrinsic properties of a surface, which implies a
series of important results. In particular, the \emph{relative curvature} of a domain on a convex surface (the ratio
of the area of the spherical image to the area of the domain) is an isometric invariant.
So the Aleksandrov's  class of convex surfaces of bounded relative curvature admits an intrinsic-geometric
definition. A.D.Aleksandrov proved that \emph{a complete convex surface of bounded relative curvature whose metric is
not everywhere Euclidean is smooth, and an incomplete surface can be non-smooth only along straight edges
with the endpoints on the boundary}. This was the first result on the dependence of the regularity of a convex surface
on the regularity of its intrinsic metric.

In contrast to the Gauss theory of surfaces, where the analytic methods dominate, in
the Aleksandrov's theory, the central role is played by the geometric methods. The main tool is the
approximation of the surface by so-called manifolds with polyhedral metric (the corresponding method in extrinsic
geometry is the approximation of a general convex surface by convex polyhedra). A simple and
natural idea of polyhedral approximation enabled one to prove the results first for polyhedra, and then, passing
to the limit, in the general case. These are the main tools of the curvature theory of general convex surfaces
which was briefly sketched above.

Let $P$ be a closed convex polyhedron. Suppose that it is cut by polygonal lines into $n$ parts,
each of which is then unfolded flat to a planar polygon $G_i, \; 1\le i \le n$. The system of
polygons $G _ {1},\dots,G_{n}$, with a given identification of their vertices and sides, is called
\emph{the net} of the polyhedron $P$. Clearly, each net satisfies the following conditions: 1) the
complex $\bigsqcup_i G_i/$identification is homeomorphic to the sphere; 2) the identified edges are
equal; 3) the sum of the angles at the identified vertices is at most $ 2\pi $.

Now suppose that we are given planar polygons $G_i$ and an identification of their sides and vertices,
which  satisfies the conditions 1),2) and 3). Then, by the Aleksandrov's ``polyhedron gluing theorem",
it is possible to glue a closed convex polyhedron (possibly bending the polygons $G_{1},\dots, G_{n}$ along
straight lines). This theorem gives the answer to the Weyl problem for polyhedral metrics. Its proof uses
rather general methods based on the topological theorem on the invariance of domain. Using the same
method A.D.Aleksandrov solved the Minkowski problem of the existence of a convex polyhedron with the prescribed areas
and directions of the faces, the  problem of the existence of a polyhedron with the prescribed curvatures and many others. This method proved to be a very effective general tool to a wide class of problems in the area.

Using the ``polyhedron gluing theorem" A.D.Aleksandrov gave a surprisingly simple solution to the Weyl problem
in the most general settings: \emph{a two-dimensional metric space of positive curvature homeomorphic to the sphere
is isomeric to a closed convex surface}. This illustrates the power of the direct geometric methods.
Alexander Danilovich used the polyhedron gluing theorem as the first step in the deep development of the whole
theory of deformations of convex surfaces. He introduced the ``gluing method" based on his ``general gluing up
theorem". Below we briefly explain the main ideas of his method.

Suppose we are given a finite number of domains in two-dimensional spaces of positive curvature.
Imagine that these domains are cut out from their spaces and some parts of their boundaries are
identified (``glued together") in such a way that the resulting space is a two-dimensional manifold $M$. If the
identified parts of the boundaries have the same lengths, then we can define in a natural way an intrinsic metric
on $M$. A.D.Aleksandrov gave necessary and sufficient conditions for $M$ to be a space of positive curvature. These
conditions are simple, natural and effective. In particular, \emph{if $M$ is homeomorphic to the sphere, then it can be
realized as a closed convex surface}. This is an extremely general theorem on the existence of a convex surface glued
from the pieces of abstractly given manifolds or from the pieces of convex surfaces. Applying the gluing method
A.D.Aleksandrov proved the following local theorem: \emph{every point of a two-dimensional space $R$ has a
neighborhood isomeric to a convex surface if and only if $R$ is a space of positive curvature}. This solves the main
problem of the intrinsic geometry of convex surfaces: the metrics of convex surfaces are characterized in a purely
intrinsic way. Using the gluing method one can prove that an ovaloid with a piece removed admits deformations by
gluing in different pieces to close the ovaloid up. One can prove, for example, that a half of an ellipsoid can be
deformed to a closed convex surface, and many other similar results. In essence, the gluing method combined with the
theorems on realizability of a metric of positive curvature by a convex surface allowed to solve, in the most
general form, all the main problems of the theory of deformation of convex surfaces \cite{1,2,3}.

It will not be an exaggeration to say that A.D.Aleksandrov has created a whole new Universe, the geometry of general
convex surfaces (a popular joke at the Conference on his 75th birthday in Novosibirsk was
``A.D.Aleksandrov is not the God, but is Godlike", which also was a reference to his large beard).
But the new Universe has to be inhabited, and this is where the difficulties and problems started.

One of the difficult open problems was the problem of congruency of closed isomeric convex surfaces and of complete
noncompact isometric convex surfaces. Another problem concerned the regularity of the convex surfaces: the isometric
immersion of an analytic metric of a positive Gaussian curvature given by the Aleksandrov's theorem produces only a
$C^{1}$-regular surface. The ``gluing" method applied to a regular convex surface of a positive Gauss curvature only
guarantees the convexity and, by the Aleksandrov's theorem on limited relative curvature, the  $C^{1}$-regularity.

The Minkowski problem asks whether there exists a closed convex surface whose Gauss curvature $K (n)$ is a given
function of the outer normal $n$. Minkowski himself proved that if the integral of $n/K(n)$ over the unit sphere is
zero, then there exists a unique (up to translation) closed convex surface with the Gauss curvature $K(n)$.
However, one can say nothing about the regularity of the surface, even when $K (n)$ is analytic. Later Levi showed
that if $K (n)$ is analytic, then the resulting surface is also analytic. In connection with these results
one may ask several questions on the regularity of the resulting surface. Namely, is it true that for a regular
function $K (n)$ the Minkowski problem has a regular solution? More precisely, if the metric is of the class $C ^{k}$,
what is the regularity class of the surface? Is it true that the convex surface (not necessarily closed) is regular
provided the function $K(n)$ is regular?

Without the answers to these questions the new theory was not complete. And the person who found the answers
was A.V.Pogorelov.

\section {Early years of A.V.Pogorelov}

A.V.Pogorelov was born on March 3, 1919, in Korocha town near Belgorod (Russia). On his farther's
``farm" there was just one cow and one horse. During the collecti\-vi\-zation they were taken from
him. Once his father came to the collective-farm stable and found his horse exhausted, dying from
thirst, while the groom was drank. Vasily Stepanovich  hit the groom, a former pauper. This
incident was reported as if a kulak has beaten a peasant, and Vasily Stepanovich was forced to
escape the town, with wife Ekaterina Ivanovna, without even taking the children. A week later
Ekaterina Ivanovna has secretly returned for the children. This is how A.V.Pogorelov came to
Kharkov, where his father became a construction worker on the building of the tractor factory.
A.V.Pogorelov told me the story of how his parents have suffered during the collectivization I have
heard from him only in 2000. In my opinion, these events had a strong influence on his life and on
the way of his public behavior. He was always very cautious in expressions and liked to quote his
mother who kept saying: silence is gold. However, he never did the things contradicting his
political views. Several times, he successfully escaped becoming a member of the Communist Party
(which was almost compulsory for a person of his scale in the USSR). As far as I know, he never
signed any letters of condemnation of dissidents, but, again, any letters in their support, as
well. Several times he was elected to the Supreme Soviet of Ukraine (although, as he said later,
against his will).

The mathematical abilities of A.V.Pogorelov became apparent already at school. His school nickname was Pascal.
He became the winner of one of the first school mathematical competitions organized by the Kharkov University, and
then of several All-Ukrainian Mathematical Olympiads. Another talent of A.V.Pogorelov was the painting. The parents
did not know, which profession to choose for him. His mother asked the son's mathematics teacher for advice.
He had a look at the paintings and said that the boy has brilliant abilities, but in the time of industrialization
the painting will not give the resource for life. This advice determined their choice. In 1937,
Alexey Vasilyevich became a student of the Department of Mathematics at the Faculty of Physics and Mathematics of
the Kharkov University.

His passion to mathematics immediately drew the attention of the teachers. Professor P.A.Solovjev gave him the book
by T.Bonnezen and V.Fenchel ``Theory of convex bodies". From that moment and for the rest of his life, geometry became
the main interest of Alexey Vasilyevich. His study was interrupted by the War.
He was conscripted and sent to study at the Air Force Zhukovski Academy. But he still thinks about geometry.
In August 1943, in a letter to Professor Ya.P.Blank he says:
"Very much I regret, that I left in Kharkov the abstracts of Bonnezen and Fenchel on the convex bodies. There are
many interesting problems in geometry ``in the large"... Do you have any interesting problem of geometry
``in the large" or of geometry in general in mind?"

After graduation from the Academy in 1945, A.V.Pogorelov starts his work as a designer engineer at the Central
Aero-Hydrodynamic Institute. But the desire to finish his university education (he finished four out of five years)
and to work in geometry brings him to Moscow University. A.V.Pogorelov asks academician I.G.Petrovsky, the head
of the Department of Mechanics and Mathematics, whether he can finish his education.
When Petrovski learnt that Alexey Vasilyevich has already graduated from the Zhukovski Academy, he decided that
there was no need in the formal completion of the university. When A.V.Pogorelov expressed his interest in geometry,
I.G.Petrovski advised him to contact V.F.Kagan. V.F.Kagan asked, what area of geometry was Alexey Vasilyevich
interested in, and the answer was: convex geometry. Kagan said that this is not his field of expertise and
suggested to contact A.D.Aleksandrov who was in Moscow at that time preparing to a mount climbing expedition
at the B.N.Delone apartment (A.D.Aleksandrov was a Master of Sports on mount climbing, and B.N.Delone was
the pioneer of Soviet mount climbing).

The first audition lasted for ten minutes. Sitting on a backpack, A.D.Aleksandrov asked Alexey
Vasilyevich the following question: \emph{is it true, that on a closed convex surface of the Gauss
curvature $K \leqslant 1$, any geodesic segment of lengths at most $\pi$ is minimizing?} It took
A.V. a year to answer this question (in affirmative) and to publish the result in 1946 in \cite{4}.
The multidimensional generalization of his theorem is a well-known theorem of Riemannian geometry,
which was proved in 1959 by W.Klingenberg: \emph{on a complete simply connected Riemannian manifold
$M^{2n}$ of sectional curvature satisfying $0 <K_{\sigma} \leqslant \lambda$, a geodesic of the
length $\leqslant \pi/\sqrt {\lambda}$ is minimizing}. In the odd-dimensional case, one needs a
two-sided bound for the curvature to obtain the same result, namely $0 < \frac {1} {4} \lambda
\leqslant K_{\sigma} \leqslant \lambda$ (and the inequality cannot be improved).

Few years ago, I asked Alexey Vasilyevich, why the Soviet mathematicians at that time showed not much interest
to the global Riemannian geometry. He answered: ``We had enough interesting problems to think about".
However, as V.A.Toponogov told me later, the first person who appreciated his comparison theorem for triangles in
a Riemannian space was  A.V.Pogorelov (in my opinion, it would be more correct to call this theorem the
Aleksandrov-Toponogov theorem, since  A.D.Aleksandrov discovered and proved it for general convex surfaces in the
three-dimensional Euclidean space).

Alexey Vasilyevich became a postgraduate-in-correspondence at Moscow State University under the supervision
of professor N.V.Efimov. Having read the manuscript of the A.D.Aleksandrov's book ``Intrinsic geometry of convex
surfaces", he starts his work in the geometry of general convex surfaces.

One of the main roles of a supervisor, in the opinion of N.V.Efimov, was to inspire a post-graduate student to
solving difficult and challenging problems. I gave numerous talks both at the N.V.Efimov's and the A.V.Pogorelov's
seminars. They were very different by style. The N.V.Efimov's seminar was long gathered, then the talk lasted for
two hours or more, and the talk was always praised very warmly, so it was almost impossible to understand the real
value of the result. A.V. always started on time, very punctually. The report lasted for at most an hour. A.V. did
not like to go through the details of the proof (probably because in many cases, after the theorem was stated,
he could prove it immediately).

In the estimation of the results he was strict and even severe. For example, in 1968, three applicants for the
Doctor degree presented their theses at the Pogorelov's seminar in Kharkov. He supported only one of them,
V.A.Toponogov, and rejected the other two, who went to Novosibirsk to A.D.Aleksandrov. All three theses were later
successfully defended.

A.V. praised rarely, but when he did -- that meant that the result was really good. He had a very fast thinking,
an enormous geometric intuition, and grasped the essence of the result very fast. Many seminar participants were
afraid to ask questions not to look foolish.

In 1947, A.V.Pogorelov defended his Candidate thesis. The main result of his thesis was the following theorem:
\emph{every general closed convex surface possesses three closed quasi-geodesics} \cite{5}. This theorem generalizes
the Lusternik~-~Shnirelman theorem on the existence of three closed geodesic on a closed regular convex surface
(a quasi-geodesic is a generalization of a geodesics; both the left and the right ``turns" of a quasi-geodesic are
nonnegative; for instance, the union of two generatrices of a round cone dividing the cone angle in two halves is
a quasi-geodesic).

After defending his Candidate thesis, A.V. discharges from the military service and moves to Kharkov (probably, this
was not an easy thing to do at that time: he was discharges by the same Order of the Defence Minister, as the son of
M.M.Litvinov, the former Soviet Minister for Foreign Affairs). In one year, he defends his Doctor Thesis
on the unique determination of a convex surface of bounded relative curvature. Soon after that, he proves the theorem
on the unique determination in the most general settings \cite{6}.

\section {Rigidity of Closed Convex Surfaces.}

In 1813, Cauchy proved the following remarkable uniqueness theorem.

\begin{named}{Cauchy Theorem}
Two closed convex polyhedra, which are equally-composed (that is, whose faces
are congruent and arranged in the same order) are congruent
\end{named}

A.D.Aleksandrov proved that it is possible, without changing the Cauchy's proof, to replace the assumption of being
equally composed with a weaker assumption of being isometric (it is easy to see that the convexity assumption
cannot be dropped: consider, for example, a cubic house with a four-chute roof and the same house with the roof
``pushed inside"). The fact that a convex surface is uniquely determined by its metric was proved for $C^{3}$-regular
closed convex surfaces of positive Gauss curvature by S.Kohn-Fossen in 1923 and for  $C^2$-surfaces
of non-negative Gauss curvature by Herglotz in 1942.

A.V.Pogorelov proved the generalization of the Cauchy Theorem to the case of arbitrary convex
surfaces:

\begin{theorem}[A.V.Pogorelov, 1949, \cite{Pog3}]
Two closed isomeric convex surfaces in the three-dimensional Euclidean space are congruent.
\end{theorem}

The incredible power of this theorem lies in the fact that it imposes no regularity assumptions on the surfaces
whatsoever. The surfaces can have edges, conic points, etc. The only extrinsic hypothesis is the convexity
(which cannot be dropped: consider the sphere and the same sphere with a small cap cut out and reflected inside).
The proof of the Cauchy Theorem in these most general settings required more than a century, and even now, after
more than half-a-century after its publication, no simpler or shorter proof is known.

A.V.Pogorelov's proof goes as follows. Suppose that there are two non-congruent isomeric closed convex surfaces
$F_0$ and $F$. Then, using the Aleksandrov's gluing theorem and his generalization of the Cauchy theorem,
it is possible to show that in an arbitrarily small neighbourhood of $F_0$, there exists a
convex surface $F_1$ isometric and non-congruent to $F_{0}$.
The next step in the proof is the ``mixing" of surfaces. The mixing of $F_0$ and $F_1$ is a family of surfaces
$F_\lambda, \; \lambda \in [0,1]$, where the surface $F _ {\lambda}$
consists of the points in the space which divide the segments, connecting the points of $F_0$ and $F_1$
corresponding by the isometry in the ratio $ \lambda: (1-\lambda)$. Then for some $\lambda$ close to $\frac12$
the surfaces $F_{\lambda} $ and $F_{1-\lambda}$ appeare to be isomeric.
A convex surface is said to be in a canonical position if it is a graph of a convex function over the $xy$-plane
and all its support planes form an angle less than $\pi/2$ with that plane. Two curves lying on two isometric
surfaces in a canonical position, are called \emph{normal equidistant curves}, if they correspond to each other by the
isometry and the corresponding points of the curves are on the same distance from the $xy$-plane. The contradiction,
which proves the theorem, is as follows: on one hand, as it can be proved, non-congruent isometric convex surfaces in
a canonical position cannot contain normal equidistant curves; on the other hand, such curves can be explicitly
produced on $F_{\lambda}$ and $F_{1-\lambda}$. This requires the theory of curves of the bounded rotation variation,
developed by A.D.Aleksandrov and V.A.Zalgaller, and also numerous geometric synthetic constructions introduced by
A.V.Pogorelov. The need for these constructions arose from the lack of analytic tools to study
conic and edge points on a convex surface. All that makes the proof extremely difficult to comprehend (in 1970,
when I was giving my first talk on a seminar in Leningrad, one of the questions was, whether I managed to
go through all the details of the proof of the Pogorelov's theorem). Another proof of this theorem follows from the
rigidity theorem of closed convex surfaces, that was also proved by A.V.Pogorelov \cite{6}. Yet another proof
can be deduced from the estimation of a deformation of a closed convex surface under a deformation of its metric
found by Yu.A.Volkov \cite{Vol}.

For solving the problem of the unique determination of a convex surface by its metric Alexey Vasilyevich was awarded
the Stalin prize of the second degree. Once, in 1951, an unexpected telegram from Kiev has arrived
saying that the Academy of Sciences nominated A.V.Pogorelov as a Correspondent Member. N.I.Ahiezer said:
``Let us pretend that we did not receive the telegram. The University itself will nominate you".

An {\it infinitesimal bending} of a surface is an infinitesimal isometric deformation. The corresponding
deformation field is called the \emph{bending field}. A bending field is called {\it trivial} if it is the
derivative at
zero of (a differentiable) rigid motion of a surface. A surface is called \emph{rigid}, if every its bending field is
trivial. W.Blaschke proved that a closed regular convex surface of a positive Gauss is rigid.

Let $F$  be a regular surface with the position vector $r=r (u, v)$ and let $\tau (u, v)$ be a smooth vector field.
The deformation $r=r (u, v) +t\tau (u, v)$ is an infinitesimal bending if and only if
$$
\langle r_u, \tau_u\rangle=0, \qquad \langle r_u,
\tau_v\rangle +\langle r_v, \tau_u\rangle=0, \qquad \langle r_v, \tau_v\rangle=0.
$$
For the $z$-component $\zeta$ of the field $\tau$ we obtain the equation
\begin{equation}\label{defo}
z _ {xx} \zeta _ {yy}-2z _ {xy} \zeta _ {xy} +z _ {yy} \zeta _ {xx} =0,
\end{equation}
which implies that the surface $z =\zeta (x, y) $ has non-positive Gauss curvature provided the curvature of $F$ is
positive.

For an infinitesimal bending of a general convex surface, equation \eqref{defo} holds almost everywhere
and we have the following lemma.

\begin{named}{Main Lemma}
If $z=z (x, y) $  is a convex surface containing no plane domains and
$\zeta (x, y)$ is the $z$-component of its infinitesimal bending field, then the surface $z=\zeta (x, y)$
has non-positive curvature, in the sense that it contains no points of the strict convexity.
\end{named}

In the regular case, this means that the Gauss curvature is non-positive. Using this fundamental fact
A.V.Pogorelov proved the following theorem.

\begin{theorem}\label{bend}
Every closed convex surface without plane domains is rigid, that is, the only possible
infinitesimal bending field are of the form $\tau=a \times r+b$, where $r$ is the position vector of the
surface, and $a, b$  are constant vectors. A closed convex surface containing plane domains is rigid outside
these domains.
\end{theorem}

A.V.Pogorelov's results on the unique determination and on the rigidity of convex surfaces formed the basis
of the geometric theory of shells. As far as I know, the physicists first disagreed with his theory, some in quite
an aggresive way. As E.P.Senkin (my PhD supervisor) was telling, when A.~V. showed them the results of the experiments
which confirmed his theory on the bending of shells they said that the shells are ``pressed by a finger".

A.V.Pogorelov's moving to Kharkov was really successful. N.I.Ahiezer drawn A.V.'s attention to the
S.N.Bernstein's papers on the Dirichlet problem for elliptic equations \cite{7}. Combining the
analytic results of S.N.Bernstein with the synthetic geometric methods A.V.Pogorelov managed to
solve the problem of the regularity of a convex surface with a regular metric of positive curvature and
of the regularity of a convex surface obtained as the solution of the Minkowski problem, with a regular positive
Gauss curvature $K(n)$.

Up to 1934, S.N.Bernstein worked in Kharkov, but then, after publishing a paper against the groundless
usage of Marxism in mathematics, he was forced to leave. After that, a public persecution of S.~N.~Bernstein began.
It worse saying that the first All-Union Mathematical Congress was held in Kharkov thanks to the fact that
S.N.Bernstein worked here.

One may regard A.V.Pogorelov as an S.N.Bernstein's successor in the field of differential equations.
Quite often, A.V. used the Bernstein's remarkable theorem which says that a non-parametric surface
of non-positive Gauss curvature defined over the whole plane and with a slower than a linear growth,
is a cylinder.

\section {Regularity of convex surfaces with a regular metric the and Min\-kowski problem.}

The regularity of a surface is equivalent to the regularity of the solutions of the Darboux equation
$$ (z _ {uu}-\Gamma ^ {1} _ {11} z_u-\Gamma ^ {2} _ {11}
z_v) (z _ {vv}-\Gamma ^ {1} _ {22} z_u-\Gamma ^ {2} _ {22} z_v)- (z _ {uv}-\Gamma ^ {1} _ {12}
z_u-\Gamma ^ {2} _ {12} z_v) ^2  =K (E-z ^ {2} _ {u}) (G-z ^ {2} _ {u}) - (F-z_uz_v) ^2,
$$
where $E, F, G$  are the coefficients of the first fundamental form, $\Gamma^{k}_{ij}$ are the
Christoffel symbols, and $z$ is the $z$-coordinate of the position vector. Its coefficients are determined
only by the metric of the surface. This nonlinear equation is the Monge-Amp\`{e}re elliptic equation
provided the Gauss curvature is positive. The question is how the regularity of solutions, for a convex surface $F$,
depends on the regularity of coefficients. A.V.Pogorelov split the solution into the three stages.
At the first stage, he considered a cap $C$, the intersection of $F$ with a closed half-space bounded by a plane $L$.
He obtained the estimates depending only on the metric for the angle between the tangent plane to
$C$ and the plane $L$, and also for the normal curvatures of $C$. The estimates for the normal curvatures at the
interior points of $C$ depend only on the metric and on the distance to $L$. To estimate the normal curvatures, he
used the method of auxiliary functions going back to S.N.Bernstein. The main difficulty, which was successfully
overcame by A.V., is to choose the correct auxiliary function for every specific problem.
These estimates correspond to the a priori estimates of the first and the second derivatives of a solution of the
Darboux equation depending on the coefficients of the equation and their derivatives.
This means that, assuming the regularity of a solution of the Darboux equation, one gives the estimates of the
derivatives of that solution depending on the coefficients and their derivatives, the distance from the boundary
and the derivatives of the solution of lower orders.
The key point is that the a priori estimates for the first and the second derivatives
are obtained from the geometric arguments (similar estimations were obtained by A.~V. in many other problems, as well).
From this point on, one obtains the a priori estimates for the higher derivatives using only analytic methods.
Let a function $z=z (x, y)$ in a domain $G$ satisfies an elliptic PDE
\begin{equation}\label{pde}
F (x, y,z, z_x, z_y, z _ {xx}, z _ {yy}, z _ {xy}) =0.
\end{equation}
Then the estimates for
the third derivatives of $z$ at a point $(x, y)$ depend on the distance from $ (x, y) $ to the boundary of
$G $, and also on
the suprema of the modules of the function $z $ and its derivatives up to the second order, the suprema of the
modules of the derivatives of the function $F $ up to the third order, and the suprema of the modules of
$(F_r)^{-1}, \, (F_t)^{-1}, \, (F_rF_t-\frac14 F^2_s)^{-1}$, where $r=z_{xx}, \, s=z_{xy}, \, t=z_{yy}$.
The estimates for the fourth and the subsequent derivatives of $z$ are based on the Schauder theory of the linear
elliptic PDE's. If $G$ is the disc $x^2+y^2\leqslant \varepsilon^2$, then the estimates on the $k$-th
derivatives of $z$ when $k \ge 3$ depend on the suprema of the module of $z$ and its derivatives up to the second
order in $G$, the suprema of the modules of $(F_r)^{-1}, \, (F_t)^{-1}, \, (F_rF_t-\frac14 F^2_s)^{-1}$,
and the suprema of the modules of the derivatives of $F$ up to order $s$, where $s=3$ for $k=3$ and
$s=k-1$ for $k> 3$. Moreover, the same values allow one to obtain the estimates for the least H\"{o}lder's
constants for the $k$-th derivatives of $z$, with an arbitrary exponent $ \alpha\ (0 <\alpha<1)$.

On the second stage, an analytic metric $g$ of positive Gauss curvature with positive geodesic
curvature of the boundary defined in a disc is realized as a convex analytic cap, which can then be
analytically extended over the boundary. To do that, A.V.Pogorelov uses the method of prolongation
over parameter,  the idea of which goes back to S.N.Bernstein, and which then was brilliantly
applied by Alexey Vasilyevich in many other problems. Following this method one has to
\begin{enumerate}[(I)]
  \item
  prove that it is possible to include the metric $g$ into a one-parameter family of analytic
metrics $g_t, \; 0\leqslant t\leqslant 1, \; g_1=g $, where the metric $g_0 $ is realized as a convex
analytic cap (in fact, $g_0$ is the metric of a spherical cap);

  \item
  prove that if the metric $g _{t_{0}}$ can be realized as a convex analytic cap, then the same is true for close
  metrics $g_t$;

  \item
  prove that if the metrics $g_{t_n} $ can be realized and $t_n \to t_0$, then the metric $g_{t_0}$ can also
  be realized.
\end{enumerate}
This will imply that the metric $g=g_1$ can be also realized as an analytic convex cap.
Problem (II) is equivalent to solving the boundary value problem for the Darboux equation
inside a unit disk with the zero boundary conditions. Let $G$ be a bounded domain with an analytic boundary
and let $\phi$ be a function analytic on $\partial G$ with respect to the arc-length parameter. By a theorem of
S.N.Bernstein, the boundary value problem for an elliptic equation $F=0$ in $G$ with $z_{|\partial G}=\phi$
has a solution if and only if the equation can be included into
a family of equations $F_t=0$ depending on a parameter $0 \leqslant t \leqslant 1 $ such that
\begin{enumerate}[(1)]
  \item $F_1=F$;
  \item the equation $F_0=0$ with the same boundary conditions has an analytic solution;
  \item for all $0 \leqslant t \leqslant 1 $, the existence of a solution $z_t$ to the boundary value problem for the
  equation $F_t=0$ with the same boundary conditions implies the uniform boundedness of $z_t$ and
  all its partial derivatives up to the second order.
\end{enumerate}
From the a priori estimates for the first and the second derivatives it follows that condition (3) of the Bernstein
theorem is satisfied, which solves (II).

Problem (III) can be solved using the a priori estimates on the derivatives up to the fourth order. This shows that
the limiting surface is $C^3$-regular. The analyticity of the limiting cap now follows from the Bernstein's theorem on
analyticity of solutions of an elliptic equation.

The third stage of the proof of the regularity of a convex surface with a regular metric goes as follows. A theorem
of A.D.Aleksandrov, the regularity of the metric and the fact that the curvature is positive imply that
any surface realizing the given metric is smooth and strictly convex. Therefore,
at any point of the surface, one can cut off a cup $F_0$ by a plane parallel to the tangent plane at that point.
Then one approximates the metric of a cap by analytic metrics in domains bounded by analytic curves of positive
geodesic curvature. These metrics can be realized as analytic caps, which converge to a limiting cap $F_1$
isomeric to $F_0$. The a priori estimates imply the regularity of the cap $F_1$. For the $C^2$ and the
$C^3$-regularity, one uses the Heintz's  a priori estimates for the first and the second derivatives of the
position vector of the surface. If, instead of the above a priori estimates for higher derivatives, one applies
the Nirenberg's theorem on the regularity of a twice differentiable solution of an elliptic equation with regular
coefficients to the Darboux equation, the results will be more precise. It follows from the Nirenberg's theorem
that a $C^2$-surface $F_1$ with a $C^n$-metric belongs to the class $C^{n-1, \alpha}, 0 <\alpha <1 $.
As $F_0$ and $F_1$ are congruent, the cap $F_0 $ is also $C^{n-1, \alpha}$-regular.

This establishes the following theorem.

\begin{theorem}[A.V.Pogorelov, \cite{Pog1,Pog2,Pog4,6}] A convex surface
with a $C^n$-regular metric, $n \geqslant 2$, of positive Gauss curvature is $C^{n-1, \alpha}$-regular,
for all $\alpha \in (0,1)$. If the metric is analytic, then the surface is also analytic.
\end{theorem}

A.V. published his first result on the regularity of convex surfaces in 1949 \cite{Pog1}. Later, in 1950,
he proved Theorem 3 for $n \ge 4$ \cite{Pog2}. The final result for $n=2,3$, was published in \cite{Pog4}.

Later, in 1953, L.Nirenberg (being familiar with the results of A.V.Pogorelov) proved the following
regularity theorem \cite{Nir1}. A $C^{m, \, \alpha}$-metric of positive Gaussian curvature, with
$m\geqslant4, \; 0 <\alpha <1$, can be realized by a $C^{m, \, \alpha}$-regular surface;
a $C^4$-metric of positive Gaussian curvature can be realized by a $C^{3, \alpha}$-regular surface, with
$0 <\alpha <1$. In the classes $C^k$, the Nirenberg's results on the regularity of the surface
are the same as the Pogorelov's ones, but in the H\"{o}lder classes they are more precise. They are
based on the a priori estimates of the H\"{o}lder class of the second derivatives of a solutions
of an elliptic equation \eqref{pde} obtained by L.Nirenberg in \cite{Nir2}.

Note that the Pogorelov's regularity theorem is fundamentally more general than the Nirenberg's one, in the
following sense.
The fact that a two-dimensional manifold with an intrinsic metric of non-negative curvature
can be locally isometrically immersed as a convex surface follows from the A.D.Aleksandrov's theorem.
Riemannian metrics of positive Gauss curvature are of that type. The Pogorelov's theorem guarantees,
that \emph{all} these surfaces are regular provided the metric is regular, while the Nirenberg's theorem says
that among them, \emph{one can find} a regular surface. This more general result was obtained with the help of
the Pogorelov's theorem on the unique determination of a convex cap, which was not used by Nirenberg.

If one considers metrics in a H\"{o}lder class, there is no loss of regularity at all:

\begin{theorem*}[I.H.Sabitov \cite{8}]
A convex surface with a $C^{n, \alpha}$-regular metric of positive curvature,
where $n\geqslant 2, \; 0 <\alpha <1$ is $C^{n, \alpha}$-regular.
\end{theorem*}

A natural question is whether the converse is true, more precisely, what is the connection between the regularity
class of a submanifold in a Riemannian space and the regularity class of the induced metric.
At first glance, the regularity of the metric should be lower. However, using the harmonic coordinates
I.H.Sabitov and S.Z.Shefel proved the following theorem.

\begin{theorem*}[\cite{9}]
Every $C ^ {k, \alpha} \ (k\geqslant 2, \ 0 <\alpha <1) $ regular
submanifold $F^l $ in a Riemannian space $M^n$ of the regularity not lower than that of $F^l$,
is a $C^{k, \alpha}$ isometrically immersed Riemannian manifold $M^l$ of the class $C^{k, \alpha}$.
\end{theorem*}

The Pogorelov's regularity theorem implied new results on the regularity of solutions of the Monge-Amp\`{e}re equation,
which became a foundation of the geometric theory of both the two-dimensional and the multidimensional theory of
the Monge-Amp\`{e}re equation (which we will discuss in Section 7).

Once I said A.V. that in my opinion, the regularity theorem is his best result, but he answered that he regards the
theorem on the unique determination as the best.

Simultaneously with the regularity theorem, in 1952 A.V. published the solution to the Minkowski problem \cite{Pog5}.

\begin{theorem}
A convex surface whose Gauss curvature is positive and is a $C^m$ function
of the outer normal $(m\geqslant 3)$, is  $C^{m+1}$-regular.
\end{theorem}

Note that this is a local theorem on a surface domain whose Gauss image is a small
disc on the unit sphere. This theorem combined with the Minkowski uniqueness theorem implies
the regularity of the solution to the Minkowski problem.
\begin{theorem}
Let $K(n)$ be a positive $C^k$ function on the unit sphere $\Omega$, $k\geqslant 3$, such that
$$
\int\nolimits_{\Omega}\frac {n \, d\omega} {K (n)} =0,
$$
where $d\omega $ is the area density on $\Omega$. Then there exists a $C^{k+1}$-regular
surface $F$ whose Gauss curvature at the point with the outer normal $n$ is $K(n)$.
The surface $F$ is unique up to a parallel translation.
\end{theorem}

In 1953, L.Nirenberg proved a similar result (using the prolongation over parameter and the a
priori estimates of Miranda): if $K \in C ^ {k,\alpha}, \; k\geqslant 2, \, 0 <\alpha <1$, then the
surface is $C^{k+1,\alpha}$-the regular. If $K \in C^2$, then the surface is $C^2$-regular
\cite{Nir1}. The Nirenberg's theorem is global, it requires the function $K$ to be defined over the
whole sphere, as he did not prove any \emph{local} theorem. Note however that the regularity is a
local property. A.V.Pogorelov told me that by the R.Courant's opinion, his results on the problem
of regularity of a surface with a regular intrinsic metric and the problem of regularity of a
solution to the Minkowski problem are more general and more natural than those of L.Nirenberg.

\section {Convex surfaces in Riemannian space.}

Perhaps the greatest achievement of A.V.Pogorelov in the area of application of the analytic methods to the theory of
convex surfaces is the following theorem.

\begin{theorem}[\cite{6,Pogor}]\label{riem}
Let $R$ be a complete three-dimensional Riemannian space whose curvature is less than some constant $C$,
and let $M$ be a Riemannian manifold homeomorphic to the sphere whose Gauss curvature is greater than $C$.

Then $M$ admits an isometric immersion in $R$. If the metrics of $R$ and $M$ are $C^n$-regular, $n\geqslant 3$,
then all such immersions are $C ^ {n-1, \alpha}$-regular, with any $\alpha \in (0,1)$.

Moreover, the isometric immersion is unique in the following sense: given any two points $x \in M$, $y \in R$, a
two-dimensional subspace $L \subset T_yR$ and a unit normal $n$ to $L$ at $y$, there exists a unique isometric
immersion $f: M \to R$ such that $f(x)=y, \; df(T_xM) = L$, and a neighborhood of $y$ on the immersed surface
$f(M)$ lies from the side of $L$ defined by $n$ (that is, the curvature vectors of all geodesic of $f(M)$ at $y$
are nonnegative multiples of $n$).
\end{theorem}

The proof of this theorem uses the prolongation over parameter and consists of three steps.
\begin{enumerate}[(I)]
  \item
First, one proves the existence of a continuous family of Riemannian manifolds $M_t, \; t \in [0,1]$, each of the
Gauss curvature greater than $C$, such that $M_1=M$ and that $M_0$ is isometrically immersible in $R$.
The manifold $M_0$ is a geodesic sphere of a small radius in $R$. Using the Aleksandrov's immersion theorem and the
Pogorelov's theorem on the regularity of a convex surface with a regular metric of the Gauss curvature greater
than $C$, one can immerse both $M_0$ and $M$ in the space of constant curvature $C$ as closed regular convex
surfaces $F_0 $ and $F$ respectively. Then it is possible to include them in a continuous family of
regular closed convex surfaces $F_t, \; t \in [0,1], \; F_1=F$, of Gauss curvature greater than $C$. Then for every
$t$, the manifold $M_t$ is $F_t$, with the induced metric.

  \item
The next step is to show that if a manifold $M_{t_0}$ is isometrically immersible in $R$, then the nearby manifolds
$M_t$ also are. First, Pogorelov considers infinitesimal bending of a regular surface in a Riemannian space. A vector
field $\xi$ on a surface $F$ in a Riemannian space $R$ is a field of infinitesimal bending if
$$
D_i\xi_j+D_j\xi_i=0, \qquad i, j=1,2,
$$
where $D_i$ is the covariant derivative in $R$, and $\xi_i$ are the covariant component the of $\xi$.
If $F$ is parameterized in such a way that its second fundamental form is $\nu ((du^1) ^2 + (du^2) ^2)$,
then the equations of the infinitesimal bending are of the form
$$
\left\{\begin {array} {l}
 {\frac {\partial \xi_1} {\partial u^1}-\frac {\partial
\xi_2} {\partial u^2} - (\tilde {\Gamma} ^ i _ {11}-\tilde {\Gamma} ^
i _ {22}) \xi_i=0,} \\
\  \\
\frac {\partial \xi_1} {\partial u^2} + \frac {\partial \xi_2} {\partial u^1}-2\tilde
{\Gamma} ^ i _ {12} \xi_i=0,
\end{array} \right.
$$
where $ \tilde {\Gamma}^i_{jk} $ are the Christoffel symbols of $F$.
The following theorem generalizes Theorem~\ref{bend} for an arbitrary ambient space.

\begin{theorem}[\cite{6,Pogor}]
Let $F$ be a convex surface homeomorphic to the sphere in a Riemannian space. Suppose $F$ has positive
extrinsic curvature. Then any field of infinitesimal bending, which vanishes at some point of $F$ together with its
covariant derivatives at that point, vanishes identically.
\end{theorem}

Let $F$ be a surface homeomorphic to the sphere, with positive extrinsic curvature in a
Riemannian space $R$. Let $F_t, \; t \in [0,1]$, be a regular deformation of $F=F_0$, and let
$ds_t^2=ds^2+td\sigma_t^2$, be the induced metric on $F_t$, where $ds^2 $ is the metric on $F$. When
$t\to 0, \; d\sigma^2_t$ tends to some limit, which is uniquely determined by the deformation. We are interested
in the inverse problem: given the limit $\lim_{t \to 0}d\sigma^2_t=d\sigma^2=\sigma _ {ij} du^idu^j$, find the
corresponding field of deformation $\xi$. The components of $\xi$ satisfy the following system of  differential
equations
$$
\left\{
\begin {array} {l}
 {\frac {\partial \xi_1} {\partial u^1}-\frac {\partial
\xi_2} {\partial u^2} - (\tilde {\Gamma} ^ i _ {11}-\tilde {\Gamma} ^ i _ {22})
\xi_i =\frac {\sigma _ {11}-\sigma _
{22}} {2},} \\
\ \ \ \ \ \\ \frac {\partial \xi_1} {\partial u^2} + \frac {\partial \xi_2} {\partial u^1}-2\tilde
{\Gamma} ^ i _ {12} \xi_i =\sigma _ {12},
\end{array}\right.
$$
which is precisely the system of differential equations for the generalized analytic functions studied by I.N.Vekua.
Using his theory Pogorelov proved that the solution $\xi$ exists and is regular. The field $\xi$ is then used in the
iterative method to find isometric immersions of metrics $M_t$ close to the immersed one $M_{t_0}$.

  \item
  The last step is to prove that if every $M_{t_n}$ is isometrically immersible, and $t_n \to t_0$, then
  $M_{t_0}$ also is.

To prove that, Pogorelov obtained the estimates on the normal curvatures of a convex surface homeomorphic to
the sphere of positive extrinsic curvature in a regular Riemannian space. These estimates depend only on the metric
of the surface and the metric of the space and, in turn, allow one to estimate the second derivatives of the
position vector of the surface. Then the estimates for the higher derivatives follow from the equation of isometric
immersion, which is elliptic, as the extrinsic curvature is positive.

Combining these three steps A.V.Pogorelov gives a solution to the generalized Weyl problem for an
isometric immersion in a Riemannian space.
\end{enumerate}

When a difficult problem is solved, then first everyone admires the solution, then gets used to it, and then, if
the theorem is not an instrumental tool, people begin forgetting it. But this is not what happened to
Theorem~\ref{riem}:
in 1997, in his talk on receiving the AMS prize for ``Pseudo-holomorphic curves in symplectic manifolds", M.Gromov
said that the idea of the proof of the existence of pseudo-holomorphic curves came to him when he was reading
the A.V.Pogorelov's proof.

As far as I know, the problem of the regularity of a convex surface with a regular metric in Riemannian
space, when the Gauss curvature of the surface is greater than the sectional curvature of the ambient Riemannian
space, is still open.

A.V.Pogorelov liked concrete problems. A.L.Verner told that when A.V. was giving a talk on the A.D.Aleksandrov's
seminar in Leningrad, he repeated several times "You Alexander Danilovich is who poses the problems, but I am who
solves them". After Pogorelov received the Lenin's prize, A.D.Aleksandrov said joking that ``we prove theorems together,
but receive prizes separately".

Alexey Vasilyevich did not have many postgraduate students. He started to work with the postgraduates students when
the Department of Geometry of the Institute of Low Temperatures was opened and the vacancies needed to be filled in.
Often, he was giving to his student a problem, the answer to which (and the method of obtaining the answer), was
already known to him. Usually, the problem concerned the improvement of some of his results. The last A.V.'s
postgraduate defended his thesis in 1970. The person who really supervised the A.V.'s postgraduates was E.P.Senkin
(he moved from Leningrad to Kharkov that time). He was a versatilely talented person who possessed a gift of praising
the students, unlike A.V.Pogorelov. Unfortunately, Eugeny Polikarpovich was ill by that time and could not help me
with the choice the problem for my thesis. In 1970, when I was the first year postgraduate, after one of the seminars,
I asked A.V., if he has a ``good" question in mind for me. He answered: ``I would be happy if you gave me the same
advice". In 1979, on the 60th anniversary of A.V., I reminded this answer to him and thanked for believing in me and
not giving me a problem for the thesis which leads to nowhere.

I never was a postgraduate student of A.V.Pogorelov, but was learning from him on his seminars what
a good problem and what a good theorem is. During more than 30 years, I gave talks at his seminar
and always expected the A.V.'s evaluation with thrill. In my 5-th year, I constructed an example of
an infinite convex polyhedron and a ray on it whose spherical image had an infinite length. There
was a conjecture that an interior point of a shortest geodesic has a neighborhood whose spherical
image has a finite length. My example provided a partial counterexample to this conjecture and I
wanted to give a talk on it on the conference on ``geometry in the large" in Petrozavodsk (Russia)
in 1969. However, the organizers rejected my application on the ground that a similar example was
constructed by V.A.Zalgaller, and they suspected me in plagiarism. Finally, I was given a time for
my report. Probably, I reported badly, so A.V. just repeated my report completely. On that
conference, I first met A.D.Aleksandrov who said to me that in his opinion, convex geometry is
already ``closed". His words made a great impression on me, so I chose a completely different area
for my postgraduate research. Actually, from that time, A.D.Aleksandrov stopped his research in
this field and started to do chronogeometry, school textbooks, ethics, philosophy, etc. It was
always easier to talk and even to debate with A.D.Alexandrov rather than with A.V.Pogorelov, he was
more liberal. At the 80th anniversary of A.V.Pogorelov I said that he was more accessible at the
conferences, but was grander in Kharkov.

\section {Surfaces of bounded extrinsic curvature}

Perhaps the most conceptual results of A.V.Pogorelov are contained in his series of papers on smooth surfaces
of bounded extrinsic curvature. In my opinion, nowadays, after half a century, these works are his most cited ones.

A.D.Aleksandrov founded the theory of the general metric manifolds, which are natural generalizations of Riemannian
manifolds. In particular, he introduced a class of two-dimensional manifolds of bounded curvature. Every metric
two-dimensional manifold which locally is a uniform limit of Riemannian manifolds whose total absolute curvatures (the
integrals of the module of the Gauss curvature) are uniformly bounded is an Alexandrov's manifold of bounded curvature.

There is a natural question, which surfaces in $\mathbb{R}^3$ carry such a metric?
A partial answer was obtained by Pogorelov, who introduced the notion of surfaces of bounded curvature.
A \emph{surface of bounded curvature} is a $C^1$-surface, the area of the spherical image of which (counting the
multiplicity of the covering) is locally bounded.

He introduced the concept of a regular point of a $C^1$- surface. A regular point can be elliptic, hyperbolic, parabolic, or a points
of inflation depending on the type of the intersection of the surface with the tangent plane. Any point on a
surface of bounded curvature can be joined to any other point in a sufficiently small neighbourhood by a rectifiable
curve lying on the surface. This defines an intrinsic metric on the underlying manifold. Pogorelov
proved that the manifold with this metric is of bounded intrinsic curvature and found connections between the
intrinsic and the extrinsic curvature of the surface. For surfaces of bounded curvature, the analogue of the Gauss
theorem holds: for every Borel set $G$, the area of the spherical image equals the \emph{intrinsic curvature}
$\omega(G)$ of the set $G$. The latter is defined by $\omega (G) = \omega ^ {+} (G)-\omega ^ {-} (G)$, where
$\omega^+(G)$ (respectively $\omega^-(G)$) is the supremum (respectively the infimum) of the sums of the
positive (respectively the negative) excesses of sets of pairwise disjoint triangles in $G$.

A very close connection was found between the intrinsic and the extrinsic geometry of a surface: a complete surface
of bounded extrinsic curvature and non-negative not identically vanishing intrinsic curvature is either a closed
convex surface or an infinite convex surface. A complete surface of bounded extrinsic curvature whose intrinsic
curvature vanishes identically is a cylinder.

The first Pogorelov's paper on surfaces of bounded extrinsic curvature was published in 1953 \cite{pog15}. On the
other hand, in 1954,
J.Nash published a paper on $C^1$ isomeric immersions, which was improved by N.Kuiper in 1955. They proved
that a two-dimensional Riemannian manifold, in rather general settings, admits a $C^1$ isometric
immersion (embedding) in $\mathbb{R}^3$. Moreover, this immersion (embedding) is, in a sense, as free as is a
topological immersion (embedding) of the underlying manifold. In particular, these results have some
``counterintuitive" corollaries: a unit sphere can be $C^1$-isometrically embedded in an arbitrarily small ball in
$\mathbb{R}^3$; there exists a closed  $C^1$-embedded locally Euclidean surface in $\mathbb{R}^3$ homeomorphic
to a torus, etc. \cite{10}. These results show that for a $C^1$-surface, even with a ``good" intrinsic metric,
there is in general no apparent connection between the intrinsic and the extrinsic geometry.
For instance, a $C^1$-regular surface whose metric is $C^2$-regular and is of positive Gauss curvature, does not
have to be convex even locally.

Perhaps the class of surfaces of bounded extrinsic curvature introduced by Pogorelov is the most natural and the widest
possible one, in which the connection between the extrinsic and the intrinsic geometry is preserved under the weakest
smoothness assumptions possible. In my opinion, this is the deepest and the most difficult series of results of
Pogorelov. The proofs are the alloy of the measure theory and brilliant synthetic geometric constructions.

\section{Multidimensional Minkowski problem and the multidimensional Monge-Amp\`{e}re equation.}

Many problems of geometry ``in the large", in the analytic form, are reduced to the existence and uniqueness problems
for certain partial differential equations, in particular, for the Monge-Amp\`{e}re equation, and conversely,
the geometric methods and results can be used to prove the existence and uniqueness of solutions for differential
equations.

As A.V.Pogorelov said about the Monge-Amp\`{e}re equation, ``This is a great equation, which I had a privilege
to study". He pioneered the study of the properties of solutions of the general multidimensional Monge-Amp\`{e}re
equation in the series of papers \cite{Pog10}-- \cite{Pog14} published in 1983-1984, and later in the monograph
\cite{12}. Earlier, he obtained the results in the case when the right-hand side is a function of the
independent variables $x_1, \ldots, x_n $, but not of the unknown function $z$ and its derivatives \cite{11}.

The Monge-Amp\`{e}re equation is a partial differential equation of the form
\begin{equation}\label{ma}
\det (z _ {ij}) =f (z_1, \dots, z_n, z, x_1, \dots, x_n), \quad f> 0,
\end{equation}
where $z_i =\frac {\partial z} {\partial x_i}, \ \ z _ {ij} = \frac {\partial^2 z} {\partial x_i \,\partial x_j}$ .
On convex solutions $z (x_1, \dots, x_n)$, this equation is of elliptic type. Rewrite equation \eqref{ma} in the form
\begin{equation}\label{matheta}
\theta (z_1, \dots,  z_n, z, x_1, \dots,  x_n) \det ( z _ {ij} ) = \varphi (x_1, \dots,  x_n).
\end{equation}
Equation \eqref{matheta} can be written in equivalent form:
\begin{equation}\label{maint}
\int\nolimits_m\theta(z_1, \dots,  z_n, z, x_1, \dots,  x_n) \det ( z _ {ij} ) dx_1 \dots  dx_n =
\int\nolimits_m\varphi (x_1, \dots,  x_n) dx_1 \dots  dx_n,
\end{equation}
where $m$ is an arbitrary Borel subset of the domain $G$ where the solution $z$ is sought. If the
solution $z (x)$ is a convex function, we can make the change of variables
$p_i =\frac {\partial z} {\partial x_i}, \ i=1, \dots, n$, on the  left-hand side. Then the resulting
equation, in contrast to equation \eqref{matheta}, makes sense for any convex but not necessarily regular
function $z(x)$. This enables one to define a \emph{generalized solution} of the  Monge-Amp\`{e}re equation as follows.
Let $z(x)$ be a convex function given in a domain $G$ and let $m \subset G$ be a Borel subset. Let $m^*$ be the set
of $p=(p_1, \dots, p_n)$ such that the hyperplane $z=p_1x_1 + \dots  +p_nx_n+c$ is a support hyperplane to the
hypersurface $z=z(x)$ at some point $(x, z(x)), \; x \in m$. Such a point $(x(p), z(p))$ is unique for almost all
$p \in m^*$. The function $z(x)$ is called a \emph{generalized solution} of equation \eqref{matheta}, if for any
Borel subset $m \subset G$,
\begin{equation*}
\int\nolimits_{m^*}\theta (p_1, \dots, p_n, z(p), x_1 (p), \dots,  x_n (p)) dp_1 \dots  dp_n =
\int\nolimits_m\varphi (x_1, \dots,  x_n) dx_1 \dots  dx_n.
\end{equation*}
The expression on the left-hand side (for an arbitrary convex function $z$) is called the \emph{conditional curvature}
of the set $m$.

The concept of a generalized solution goes back to A.D.Aleksandrov. One of the main problems is to prove the existence
of a generalized solution of \eqref{matheta} under certain natural assumptions about the functions $ \varphi, \theta $.
In the two-dimensional case with $\theta=1$, the existence was proved by A.V.Pogorelov, and in the general case, by
A.D.Aleksandrov. Then Pogorelov proved the existence of a solution of the Dirichlet problem and the maximum principle
for generalized solutions of the Monge-Amp\`{e}re equation with $\theta_z \leqslant 0$, which implies the uniqueness
of a solution of the Dirichlet problem. He also considered similar problems for the Monge-Amp\`{e}re equation on the
sphere.

The first step in the proof of the existence of a generalized solution of the Monge-Amp\`{e}re
equation is to show that there exists a convex polyhedron whose vertices project to the given
points $B_i$ in the $x$-space, with the given conditional curvatures $\mu_i > 0$. Next, by passing
to the limit, one proves that given a convex domain $G$ and a bounded measure $\mu$ on the Borel
subsets of $G$, there exists a convex hypersurface $z=z(x), \; x \in G$, such that for every Borel
subset $m \in G$, the conditional curvature of $m$ is $\mu(m)$. For the Monge-Amp\`{e}re equation,
$\mu (m) = \int\nolimits_m \varphi (x) dx$. If the function $\theta (p, z, x)$, which defines the
conditional curvature, strictly decreases in $z$, then the convex hypersurface $z = z(x)$ is
uniquely determined by its boundary values and the measure $\mu$. This implies the existence and
the uniqueness of a generalized solution of the Monge-Amp\`{e}re equation. The proof of the
regularity of the solution assuming the functions $\theta$ and $\varphi$ to be regular, can be
reduced to the proof of the regularity of a convex hypersurface with the given conditional
curvature. This can be done using the prolongation over parameter. The most difficult part of the
proof is finding the a priori estimates for the solution and its derivatives up to the third order
(the a priori estimates for the higher order derivatives can be then obtained from equation
\eqref{matheta}). A.V.Pogorelov proved the following theorem.

\begin{theorem}
A generalized solution of the Monge-Amp\`{e}re equation \eqref{matheta}, with $\theta$ and $\varphi$ being regular
positive functions and $\theta_z\leqslant 0$, is regular in a neighborhood of every point of the strict convexity.
If $\theta, \varphi\in C^k(G), \; k\geqslant 3$, then the solution is $C^{k+1, \alpha}$-regular for all
$\alpha \in (0, 1)$.
\end{theorem}

The key role in this theorem (as in many others) is played by the a priori estimates of a solution of an elliptic
equation, together with its derivatives (up to the third order, for the Monge-Amp\`{e}re equation). These estimates
do not directly follow from the equation. They are needed to guarantee the $C^{2, \alpha}$-regularity of the
limiting solution of a sequence of regular solutions. Then, if the coefficients are regular, one can apply the
standard tools of the theory of elliptic equations.

Pogorelov's a priori estimates for the third derivatives of a solution $z=z(x)$ of \eqref{ma} are based on the
idea of E.Calabi's \cite{Calabi}, who considered a Riemannian metric $ds^2=z_{ij} dx^i dx^j$, computed the Laplacian
of the scalar curvature of it and obtained the estimates for the third derivatives of $z$ in the case
$f=\mathrm{const}> 0$.

\bigskip

Earlier, Pogorelov solved the multidimensional Minkowski problem. In 1968 -- 1971, he published a
series of papers in ``DAN", the Doklady Akademii Nauk (Proceedings of the USSR Academy of Sciences)
where he found a priori estimates for the second and the third derivatives, and proved regularity
of a solution of the Minkowski problem in the multidimensional case using the prolongation over
parameter.

\begin{theorem}[\cite{11,Pog7}]
Let $K(n)$ be a positive $C^k$ function, $k\geqslant 3$, on the unit sphere $S^{m-1} \subset R^m$ satisfying
$$
\int\nolimits_{S^{m-1}} \frac {n} {K (n)} d\omega=0.
$$
Then there exists a (unique up to parallel translation) $C^{k+1, \alpha}$-regular convex hypersurface with the
Gauss curvature $K(n)$ at the point with the outer normal $n$, for every $0 <\alpha <1$. If $K (n) $ is analytic, then
the corresponding hypersurface is also analytic.
\end{theorem}

Using this theorem Pogorelov proved the regularity of the solutions of equation \eqref{matheta} with $\theta=1$
and the regularity of generalized solutions of the Dirihlet problem \cite{Pog8,Pog9,11}. One of the implications
is the following theorem: a unique convex solution $z=x (x_1, \ldots, x_n) $ of the equation
$\det (z_{ij}) =\mathrm{const}> 0$ defined over the whole space $x_1, \ldots, x_n$ is a quadratic polynomial
\cite{11}.

Note that A.V.Pogorelov did not publish long papers from the middle of 50-th. Usually he published a brief note in
DAN, and later, a separate small book. For instance, the book ``Multidimensional Minkowski Problem" was published
only in 1975.

In 1976 -- 1977,  S.Yu.Cheng and S.T.Yau published the papers \cite{Yau1,Yau2}. They found small inaccuracies and
incompleteness in the Pogorelov's brief notes in DAN (avoiding technicalities was a matter of the journal
style; later all the details were given in \cite{11}) and declared that Pogorelov had no complete proof. Then they
gave their own proof of regularity of a solution of the multidimensional Minkowski problem and the Monge-Amp\`{e}re
equation
$$
\det (z _ {ij}) =F (x_1, \ldots x_n, z)> 0.
$$
The proof heavily relies on the Pogorelov's a priori estimates on the second and the third derivatives, and on the
other results, in particular, the Aleksandrov's and Pogorelov's theorems on the generalized solutions of the
Monge-Amp\`{e}re equation. The results of these papers can be viewed as a restatement of the earlier
results of Pogorelov, but by no means as the new results.

The English translation of ``Multidimensional Minkowski problem" appeared in 1978 with a very nice foreword by
L.Nirenberg. However, nowhere in the translation it is mentioned \emph{when the Russian original was published}.
This is an (unfortunate) reason why often the authors first refer to \cite{Yau1} and then to the English translation
of the Pogorelov's book.

\bigskip

Let $M$ be a compact K\"ahler manifold with the K\"ahler metric $ds^2=g_{j\bar {k}} dz^j d\bar{z}^k$ and the
K\"ahler form $\omega =\frac {i} {2} g _ {j\bar {k}} dz^j\wedge d\bar {z} ^k$. In 1954, E.Calabi conjectured that
for a given $(1,1)$ form $\sigma =\frac {i} {2\pi} \tilde {R}_{j\bar {k}} dz^j\wedge d\bar {z} ^k$
representing the first Chern class of $M$, there exists a K\"ahler metric on $M$ with the Ricci tensor
$\tilde {R} _ {j\bar {k}}$ whose K\"ahler form belongs to the same cohomology class as $\omega$. To solve the
Calabi conjecture, one needs to solve the complex Monge-Amp\`{e}re equation
$$
\det\Bigl(g _ {i\bar {k}} + \dfrac {\partial^2\varphi} {\partial z^j\partial\bar {z} ^k} \Bigr)
=\det (g _ {j\bar {k}}) e ^ F, \quad c=0, 1,
$$
for a real function $ \varphi $, where $F $ is a given function satisfying $\int_M e^F dz
=\mathrm{Vol} \, M$ in the case $c=0$.

S.T.Yau solved the Calabi conjecture using the prolongation over parameter. A substantial part of
the proof is finding the a priori estimates for the second and the third derivatives. In his paper
\cite{Yau3}, Yau writes that the Pogorelov's estimates for the real Monge-Amp\`{e}re equation were
the basis for his estimates in the complex case but give no references on Pogorelov's articles.

However, in the recent paper ``Perspectives on Geometric Analysis" (arXiv: math.DG/0602363), by an inaccurate citing, Yau again lessens the Pogorelov's role in the solving of two fundamental problems:
the regularity of a convex surface with a regular metric and the Minkowski problem. Concerning the first one, he
refers only to the 1961 paper in DAN, but not to the 1949 paper \cite{Pog1} or a long paper in the
``Notes of the Kharkov Mathematical Society" published in 1950. Also, he refers to the 1953 Nirenberg's paper, which
was published later than the original paper of Pogorelov. As to the two-dimensional Minkowski problem, he refers
to the 1953 Pogorelov's paper which has nothing to do with the subject (it follows even from the title). Meanwhile,
the Pogorelov's solution of the Minkowski problem was published in 1952, \cite{Pog5}, which is a year earlier than the
Nirenberg's paper cited by Yau. As to the multidimensional Minkowski problem, only the 1976 Cheng-Yau's paper is
cited. There are no references to the Pogorelov's papers and books on the subject (even those translated into English)
whatsoever. Another example is the book  \cite{Ger}. On page 256 the author says:
``The(Minkowski) Problem has been partially solved by Minkowski, Aleksandrov, Lewy, Nirenberg, and Pogorelov" then
referring to the 1952 paper \cite{Pog5}, but without mentioning the papers \cite{Pog8, Pog9} and the 1975 book
\cite{11} (English translation 1978).

The Pogorelov's results on the multidimensional Monge-Amp\`{e}re equation were generalized in various directions,
such as improving the regularity of the solution \cite{Caf}, proving the regularity up to the boundary \cite{CNS},
studying the degenerate Monge-Amp\`{e}re equations \cite{GTW}, and applying the results and the methods to other
classes of completely nonlinear second order equations \cite{SUW}. In particular, it is proved in \cite{Caf} that
a convex solution $z=z(x)$ of the equation
\begin{equation*}
\det (z _ {ij}) =f (x_1, \ldots, x_n)
\end{equation*}
in a convex domain $\Omega$, with $f\in C ^\alpha(\Omega)$, belongs to $C^{2, \alpha}(\Omega)$ (which is the best
possible smoothness).

Note that the multidimensional Monge-Amp\`{e}re equation appears also in statistical mechanics, meteorology,
financial mathematics and in other areas \cite{Caf2}.

Curiously enough, none of the Pogorelov's students in Kharkov worked in differential equations, but his methods
and results were activley developed by the Leningrad mathematical school. It seems that the mathematical influence
of A.V.Pogorelov was proportional to the distance from Kharkov.

\begin{center}
* \hspace{1em} * \hspace{1em} *
\end{center}

The mathematical legacy of A.V.Pogorelov is enormous. The most influential results not mentioned in the previous
sections include a complete solution of the Fourth Hilbert Problem and obtaining the necessary and sufficient
conditions for a $G$-space to be Finsler \cite{13,14}.

Until 1970, A.V.Pogorelov lectured at Kharkov University. Based on this lecture notes, he published a series of
brilliant textbooks on analytic and differential geometry and the foundation of geometry. Sometimes, during
routine lectures, he was thinking about his research. Anecdote says that on one of such lectures reflecting
on something completely different he started improvising and became lost. Then he opened the textbook with the words:
``What does the author say on the topic? Oh, yes, it is obvious  \dots ". In contrast, when lecturing on a topic
interesting to him, A.V.Pogorelov was very enthusiastic and inspired (I remember one of his topology courses for
the 4th year students). But perhaps the best of his lecturing brilliance was seen when he was presenting his own
results. His talks were real fine art performances. In his opinion, one of the most valuable qualities of a
mathematical result is its beauty and naturalness. That is why he usually omitted technicalities, and for the
sake of simplicity and beauty was ready to sacrifice the generality.

For many years, A.V.Pogorelov was the editor-in-chief of the ``Ukrainskij Geometricheskij Sbornik"
(Ukrainian Geometrical, a geometry
journal published annually in Kharkov University. He was very jealous about the publications in it by the ``local"
mathematicians. I remember, once I submitted a paper to a different journal, but has not submitted one to the UGS. A.V.
became very disappointed at me and said: ``Wise people say, the one who wants to become famous, must become
famous on his own place".

A.V.Pogorelov was the author of one of the most popular school textbooks in geometry. This began as follows. He was
a member of the commission on the school education whose head was A.N.Kolmogorov. A.V. disagreed with the textbook
written
by A.N.Kolmogorov and his coauthors and wrote his own manual for teachers on elementary geometry, in which he
built the whole school geometry course starting with a set of natural and intuitive axioms. The manual was
published in 1969 and formed a basis for his school textbook. A.V. used to say: ``My textbook is the Kiselyov's
improved textbook" (``Elementary geometry" by A.P.Kiselyov is probably the most well-known Russian-language school
geometry textbook; it was first published in 1892, with the last edition in 2002; many generations of students
studied the Kiselyov's ``Geometry"). The first version of the A.V.Pogorelov's textbook sparked sharp criticism from
A.D.Aleksandrov whom Pogorelov deeply respected. This criticism was based on implementing the axiomatic approach as
early as in year six at school: ``What is the point to prove `obvious' statements (from the student's point of
view)?". After reworking of the textbook, these disagreements were resolved, and they remained in strong friendship
till the last days of A.D.Aleksandrov.

Alexey Vasilyevich was a person of the highest decency. When a five year contract with the ``Prosvescheniye" Publisher
was coming to an end, another publisher offered a very tempting contract to him. He refused on the unique ground
that it will be unfair to the editor of the textbook. It should be noted that the money for the school textbook
republishing were the main source of his living in the middle of the 90-th.

A.V.Pogorelov told me that I.G.Petrovsky invited him to the Moscow University, I.M.Vinogradov invited to Moscow
Mathematical Institute, A.D.Aleksandrov invited to Leningrad several times. He even spent one year (1955-1956) in
Leningrad, but then returned to Kharkov. He preferred to stay in Kharkov, far from the fuss and noise of the capitals.
In Kharkov he proved his theorems, and to Moscow and Leningrad he went to shine.

A.V.Pogorelov is a remarkable example of the mathematical longevity. I remember, in 1992, on the A.D.Aleksandrov's
80-th anniversary, I asked M.Berger a question on geometry. He answered that he is too old for geometry. He was 65 at
that time. However, at this very age A.V.Pogorelov received his final results on the multidimensional
Monge-Amp\`{e}re equation. Only in 1995, in the age of 76, he said ``At my age, it is already not necessary to do
mathematics".

A.V.Pogorelov combined in himself a diligence of the peasant and a mathematical brilliancy. He simply could not live
without working. He inherited this from his parents, Vasily Stepanovich and Ekaterina Ivanovna. On the 50-th
anniversary of A.V.Pogorelov, N.I.Akhiezer bowed to his parents thanking them for their son.

A.V.Pogorelov was a very handsome man. He liked to be photographed. He often behave artistically and had a good
sense of humor. Once I visited A.V. and left a bag at his place. When i returned to pick it up,
A.V. said with a smile: ``No need to apologize. If you forgot, then you probably were thinking about something".
I remember as in 1982, on A.D.Aleksandrov's 70-th anniversary, somebody asked his employee A.I.Medjanik to sing
(he has a beautiful voice). After the song, A.V. said: ``Have you heard him singing? Now imagine how the boss sings!"
R.J.Barri, N.V.Efimov's wife, told me that when Pogorelov was a Efimov's post-graduate student, Efimov never invited
anyone to his place on the day when they had consultations (on Thursday). This rule was broken only once, when
V.A.Rokhlin was leaving Moscow. On that party, Alexey Vasilyevich sang the Ukrainian songs.

A.V.Pogorelov was a modest person, despite of all his titles. In 1972, when the Kharkov geometers flied through
Moscow to Samarkand to the All-Union geometry conference, the flight was delayed in Moscow and we had to
spend a night in the waiting hall. Being a Member of the Supreme Soviet of Ukraine, Pogorelov could go to a
VIP-hall, but he has chosen to stay with us.

In 2000, at the age of 81, A.V.Pogorelov moved to Moscow. He lived in Novokosino (one of the outer suburbs of Moscow),
and when going to the Institute of Mathematics, used only the public transport, with several changes, instead of
calling a car from the Academy of Sciences.
In Moscow he continued to work, to think on the geometry problems and not only on them. He even brought from Kharkov a
drawing board, on which he projected an electric generator based on the superconductivity.

Alexey Vasilyevich Pogorelov was a person blessed by an incredible natural talent combined with a constant tireless
labor.

\bigskip

\textbf{Acknowledgement.} I cordially thank V.A.Marchenko for inspiring me on this paper and for
many critical remarks and J.G.Reshetnjak for thoroughly reading the manuscript and making numerous
useful corrections. I am also thankful to A.L.Yampolsky and to Yu.Nikolayevsky for the English
translation.

\end{document}